\documentclass[11pt,a4paper]{scrartcl}
\usepackage{amsmath}
\usepackage{amssymb}
\usepackage{amsthm}
\usepackage[english]{babel}
\usepackage{verbatim}
\usepackage[shortlabels]{enumitem}
\usepackage{tikz-cd}
\usepackage{etoolbox} 
\usepackage[all]{xy}
\usepackage{hyperref}

\newtheoremstyle{mio}%
	{}{} 
	{\itshape}{} 
	{\bfseries}{.}{ } 
	{#1 #2\thmnote{~\mdseries(#3)}} 
\theoremstyle{mio}
\newtheorem{teor}{Theorem}[section]
\newtheorem{cor}[teor]{Corollary}
\newtheorem{prop}[teor]{Proposition}
\newtheorem{lemma}[teor]{Lemma}

\theoremstyle{definition}
\newtheorem{Def}[teor]{Definition}

\newcommand{\Br}{\mathrm{Br}}

\newcommand{\ins}[1]{\mathbb{#1}}
\newcommand{\insN}{\ins{N}}

\newcommand{\insQ}{\ins{Q}}
\newcommand{\insR}{\ins{R}}
\newcommand{\inN}{\in\insN}

\newcommand{\inR}{\in\insR}
\newcommand{\N}{\mathbb{N}}

\newcommand{\Zar}{\mathrm{Zar}}

\newcommand{\Q}{\mathbb{Q}}

\newcommand{\insiemeVE}{\mathcal{V}}

\newcommand{\insiemeVEdiv}{\mathcal{V}_\mathrm{div}}
\newcommand{\insiemeVEstaz}{\mathcal{V}_\mathrm{stat}}
\newcommand{\cons}{\mathrm{cons}}

\newcommand{\limiti}{\mathcal{L}}

\newcommand{\ball}{\mathcal{B}}
\newcommand{\sq}{\delta}
\newcommand{\corona}{\mathcal{C}}

\newcommand{\R}{\mathbb{R}}

\newcommand{\Spec}{\mathrm{Spec}}
\newcommand{\Max}{\mathrm{Max}}

\newcommand{\limtop}[3]{(#1,#2]_{#3}}
\newcommand{\limtopstd}{\limtop{-\infty}{\delta(\beta,K)}{\insQ\Gamma_v}}

\title{Metrizability of spaces of valuation domains associated to pseudo-convergent sequences}

\author{
G.\ Peruginelli\footnote{Dipartimento di Matematica "Tullio Levi-Civita", University of Padova, Via Trieste 63,
35121 Padova, Italy. E-mail: gperugin@math.unipd.it}\\
D. Spirito\footnote{Dipartimento di Matematica "Tullio Levi-Civita", University of Padova, Via Trieste 63,
35121 Padova, Italy. E-mail: spirito@math.unipd.it}
}

\begin{document}
\maketitle

\begin{abstract}
\noindent Let $V$ be a valuation domain of rank one with quotient field $K$. We study the set of extensions of $V$ to the field of rational functions $K(X)$ induced by pseudo-convergent sequences of $K$ from a topological point of view, endowing this set either with the Zariski or with the constructible topology. In particular, we consider the two subspaces induced by sequences with a prescribed breadth or with a prescribed pseudo-limit. We give some necessary conditions for the Zariski space to be metrizable (under the constructible topology) in terms of the value group and the residue field of $V$. 
 \\


\noindent MSC Primary 12J20, 13A18, 13F30.
\end{abstract}

\section{Introduction}

Let $D$ be an integral domain with quotient field $K$, and let $L$ be a field extension of $K$. The \emph{Zariski space} $\Zar(L|D)$ of $L$ over $D$ is the set of all valuation domains containing $D$ and having $L$ as quotient field. This set was originally studied by Zariski during its study of the problem of resolution of singularities \cite{zariski_sing,zariski_comp}; to this end, he introduced a topology (later called the \emph{Zariski topology}) that makes $\Zar(L|D)$ into a compact space that is not Hausdorff \cite[Chapter VI, Theorem 40]{ZS2}.

A second topology that can be considered on the Zariski space is the \emph{constructible topology} (or \emph{patch topology}), that can be constructed from the Zariski topology in the same way as it is constructed on the spectrum of a ring. The Zariski space $\Zar(L|D)$ endowed with the constructible topology, which we denote by $\Zar(L|D)^\cons$, is more well-behaved than the starting space $\Zar(L|D)$ with the Zariski topology, since beyond being compact it is also Hausdorff; furthermore, it keeps its link with the spectra of rings, in the sense that there is a ring $A$ such that $\Spec(A)$ is homeomorphic to $\Zar(L|D)^\cons$ \cite{hochster_spectral}.

Suppose now that $D=V$ is a valuation domain. In this case, the study of $\Zar(L|V)$ often concentrates on the subset of the \emph{extensions} of $V$ to $L$, i.e., to the valuation domains $W\in\Zar(L|V)$ such that $W\cap K=V$. When $L=K(X)$ is the field of rational functions over $K$, there are several ways to construct extensions of $V$ to $K(X)$, among which we can cite key polynomials \cite{maclane,vaquie}, monomial valuations, and minimal pairs \cite{APZ1,APZ2}. Another approach is by means of \emph{pseudo-monotone sequences} and, in particular, \emph{pseudo-convergent sequences}: the latter are a generalization of the concept of Cauchy sequences that were introduced by Ostrowski \cite{Ostr} and later used by Kaplansky to study immediate extensions and maximal valued fields \cite{Kap}. Pseudo-monotone sequences were introduced by Chabert in \cite{ChabPolCloVal} to describe the polynomial closure of subsets of rank one valuation domains. In particular, Ostrowski introduced pseudo-convergent sequences in order to describe all rank one extensions of a rank one valuation domain when the quotient field $K$ of $V$ is algebraically closed (\emph{Ostrowski's Fundamentalsatz}, see \cite[\S 11, IX, p. 378]{Ostr}); recently, the authors used pseudo-monotone sequences  to extend Ostrowski's result to arbitrary rank when the completion $\widehat{K}$ of $K$ with respect to the $v$-adic topology is algebraically closed \cite[Theorem 6.2]{PS2}.

Motivated by these results, in this paper we are interested in the subspace $\insiemeVE$ of $\Zar(K(X)|V)$ containing the extensions of $V$ defined by pseudo-convergent sequences, under the hypothesis that $V$ has rank $1$ (see \S 2 for the definition of this kind of extensions). The study of $\insiemeVE$ was started in \cite{PS1}, where it was shown that $\insiemeVE$ is always a regular space (even under the Zariski topology) \cite[Theorem 6.15]{PS1} and that the Zariski and the constructible topology agree on $\insiemeVE$ if and only if the residue field of $V$ is finite \cite[Proposition 6.11]{PS1}. We continue the study of this space by concentrating on the problem of metrizability: more precisely, we are interested on conditions under which $\insiemeVE$ and some distinguished subsets of $\insiemeVE$ are metrizable. More generally, we look for conditions under which the whole Zariski space (endowed with the constructible topology) is metrizable. To do so, we consider two partitions of $\insiemeVE$.

In Section \ref{sect:radfix}, we study the spaces $\insiemeVE(\bullet,\delta)\subset\insiemeVE$ consisting of those extensions of $V$ induced by pseudo-convergent sequences having the same (fixed) breadth $\delta\inR\cup\{\infty\}$ (see Section \ref{sect:background} for the definition); this can be seen as a generalization of the study of valuation domains associated to elements of the completion of $K$ tackled in \cite{PerTransc}, which in our notation reduces to the special case $\delta=\infty$. In particular, we show that $\insiemeVE(\bullet,\delta)$ can be seen as a complete ultrametric space under a very natural distance function (Theorem \ref{teor:Sigmadelta}) which induces both the Zariski and the constructible topology (that in particular coincide, see  Proposition \ref{prop:compl}); however, these distances (as $\delta$ ranges in $\R\cup\{\infty\}$), cannot be unified into a metric encompassing all of $\insiemeVE$  (Proposition \ref{prop:restriz-dsq}).

In Section \ref{sect:bfix}, we study the spaces $\insiemeVE(\beta,\bullet)\subset\insiemeVE$ consisting of those extensions of $V$ induced by pseudo-convergent sequences having a (fixed) pseudo-limit $\beta\in\overline{K}$ (with respect to some prescribed extension of $V$ to $\overline{K}$). We show that these spaces are closed, with respect to the Zariski topology (Proposition \ref{prop:betachiuso}), and that the constructible and the Zariski topology agree on each $\insiemeVE(\beta,\bullet)$ (Proposition \ref{prop:betafix-constr}); furthermore, we represent $\insiemeVE(\beta,\bullet)$ through a variant of the upper limit topology (Theorem \ref{teor:betafix}), and we show that it is metrizable if and only if the value group of $V$ is countable (Proposition \ref{prop:XGamma-metriz}). As a consequence, we get that, when the value group of $V$ is not countable, the space $\Zar(K(X)|V)^\cons$ is not metrizable (Corollary \ref{cor:Gammav-uncont-metriz}).

In Section \ref{sect:beyond}, we look at the same partitions, but on the sets $\insiemeVEdiv$ and $\insiemeVEstaz$ of extensions induced, respectively, by pseudo-divergent and pseudo-stationary sequences (the other type of pseudo-monotone sequences beyond the pseudo-convergent ones, see \cite{ChabPolCloVal, PerPrufer, PS2}). Using a quotient onto the space $\Zar(k(t)|k)$ (where $k$ is the residue field of $V$) we first show that $\Zar(K(X)|V)^\cons$ is not metrizable if $k$ is uncountable (Proposition \ref{prop:uncount-resfield}); then, with a similar method, we show that $\insiemeVEdiv(\bullet,\delta)$ is not Hausdorff (with respect to the Zariski topology) when $\delta$ belongs to the value group of $V$ (Proposition \ref{prop:VEdiv-delta}). On the other hand, we show that fixing a pseudo-limit (i.e., considering $\insiemeVEdiv(\beta,\bullet)$) we get a space homeomorphic to $\insiemeVE(\beta,\bullet)$ (Proposition \ref{prop:VEdiv-beta}). For pseudo-stationary sequences, we show that both partitions $\insiemeVEstaz(\bullet,\delta)$ and $\insiemeVEstaz(\beta,\bullet)$ give rise to discrete spaces (Proposition \ref{prop:VEstaz}).

\section{Background and notation}\label{sect:background}
Let $D$ be an integral domain and $L$ be a field containing $D$ (not necessarily the quotient field of $D$). The \emph{Zariski space of $D$ in $L$}, denoted by $\Zar(L|D)$, is the set of valuation domains of $L$ containing $D$ endowed with the so-called \emph{Zariski topology}, i.e., with the topology generated by the subbasic open sets
\begin{equation*}
B(\phi)=\{W\in\Zar(L|D)\mid \phi\in W\},
\end{equation*}
where $\phi\in L$. Under this topology, $\Zar(L|D)$ is a compact space \cite[Chapter VI, Theorem 40]{ZS2}, but it is usually not Hausdorff nor $T_1$ (indeed, $\Zar(L|D)$ is a $T_1$ space if and only if $D$ is a field and $L$ is an algebraic extension of $D$). The \emph{constructible topology} on $\Zar(L|D)$ is the coarsest topology such that the subsets $B(\phi_1,\ldots,\phi_k)=B(\phi_1)\cap\cdots\cap B(\phi_n)$ are both open and closed. The constructible topology is finer than the Zariski topology, but $\Zar(L|D)^\cons$ (i.e., $\Zar(L|D)$ endowed with the constructible topology) is always compact and Hausdorff \cite[Theorem 1]{hochster_spectral}.

\medskip

From now on, and throughout the article, we assume that $V$ is a valuation domain of rank one; we denote by $K$ its quotient field, by $M$ its maximal ideal and by $v$ the valuation associated to $V$. Its value group is denoted by $\Gamma_v$. 

If $L$ is a field extension of $K$, a  valuation domain $W$ of $L$ \emph{lies over} $V$ if $W\cap K=V$; we also say that $W$ is an \emph{extension} of $V$ to $L$. In this case, the residue field of $W$ is naturally an extension of the residue field of $V$ and similarly the value group of $W$ is an extension of the value group of $V$.

We denote by $\widehat{K}$ and $\widehat{V}$ the completion of $K$ and $V$, respectively, with respect to the topology induced by the valuation $v$. We still denote by $v$ the unique extension of $v$ to $\widehat K$ (whose valuation domain is precisely $\widehat{V}$). We denote by $\overline{K}$ a fixed algebraic closure of $K$. 

Since $V$ has rank one, we can consider $\Gamma_v$ as a subgroup of $\insR$. If $u$ is an extension of $v$ to $\overline{K}$, then the value group of $u$ is $\insQ\Gamma_v=\{q\gamma\mid q\in\insQ,\gamma\in\Gamma_v\}$. 

The valuation $v$ induces an ultrametric distance $d$ on $K$, defined by
\begin{equation*}
d(x,y)=e^{-v(x-y)}.
\end{equation*}
In this metric, $V$ is the closed ball of center 0 and radius 1. Given $s\in K$ and $\gamma\in\Gamma_v$, the closed  ball of center $s$ and radius $r=e^{-\gamma}$ is:
$$\{x\in K \mid d(x,s)\leq r\}=\{x\in K \mid v(x-s)\geq \gamma\}.$$

\medskip

The basic objects of study of this paper are pseudo-convergent sequences, introduced by Ostrowski in \cite{Ostr} and used by Kaplansky in \cite{Kap} to describe immediate extensions of valued fields. Related concepts are \emph{pseudo-stationary} and \emph{pseudo-divergent} sequences introduced in \cite{ChabPolCloVal}, which we will define and use in Section \ref{sect:beyond}.
\begin{Def}
Let $E=\{s_n\}_{n\in\N}$ be a sequence in $K$. We say that $E$ is a \emph{pseudo-convergent} sequence if $v(s_{n+1}-s_n)<v(s_{n+2}-s_{n+1})$ for all $n\in\N$.
\end{Def}
In particular, if $E=\{s_n\}_{n\in\N}$ is a pseudo-convergent sequence and $n\geq 1$, then $v(s_{n+k}-s_n)=v(s_{n+1}-s_n)$ for all $k\geq 1$. We shall usually denote this quantity by $\delta_n$; following \cite[p. 327]{Warner} we call the sequence $\{\delta_n\}_{n\inN}$ the \emph{gauge} of $E$. We call the quantity
\begin{equation*}
\delta_E=\lim_{n\to\infty}v(s_{n+1}-s_n)=\lim_{n\to\infty}\delta_n
\end{equation*}
the \emph{breadth} of $E$. The breadth $\delta_E$ is an element of $\insR\cup\{\infty\}$, and it may not lie in $\Gamma_v$.

\begin{Def}
The \emph{breadth ideal} of $E$ is
\begin{equation*}
\Br(E)=\{b\in K \mid v(b)>v(s_{n+1}-s_n),\forall n\in\N\}=\{b\in K\mid v(b)\geq\delta_E\}.
\end{equation*}
\end{Def}
In general, $\Br(E)$ is a fractional ideal of $V$ and may not be contained in $V$. If $\delta=+\infty$, then $\Br(E)$ is just the zero ideal and $E$ is a Cauchy sequence in $K$. If $V$ is a discrete valuation ring, then every pseudo-convergent sequence is actually a Cauchy sequence.

The following definition has been introduced in \cite{Kap}, even though an equivalent concept already appears in \cite[p. 375]{Ostr} (see \cite[X, p. 381]{Ostr} for the equivalence).
\begin{Def}
An element $\alpha\in K$ is a \emph{pseudo-limit} of $E$ if $v(\alpha-s_n)<v(\alpha-s_{n+1})$ for all $n\in\N$, or, equivalently, if $v(\alpha-s_n)=\delta_n$ for all $n\in\N$. We denote the set of pseudo-limits of $E$ by $\limiti_E$, or $\limiti_E^v$ if we need to emphasize the valuation.
\end{Def}
If $\Br(E)$ is the zero ideal then $E$  is a Cauchy sequence in $K$ and converges to an  element of $\widehat{K}$, which is the unique pseudo-limit of $E$. Kaplansky proved the following more general result.
\begin{lemma}\label{lemma:kaplansky}
\cite[Lemma 3]{Kap} Let $E\subset K$ be a pseudo-convergent sequence. If $\alpha\in K$ is a pseudo-limit of $E$, then the set of pseudo-limits of $E$ in $K$ is equal to $\alpha+\Br(E)$. 
\end{lemma}
Lemma \ref{lemma:kaplansky} can also be phrased in a geometric way: if $\alpha\in\limiti_E$, then $\limiti_E$ is the closed ball of center $\alpha$ and radius $e^{-\delta_E}$.

The following concepts have been given by Kaplansky in \cite{Kap} in order to study the different kinds of immediate extensions of a valued field $K$, i.e., extensions $V\subseteq W$ of valuation rings where neither the residue field nor the value group change. 

\begin{Def}
Let $E$ be a pseudo-convergent sequence. We say that $E$ is of \emph{transcendental type} if, for every $f\in K[X]$, the value  $v(f(s_n))$ eventually stabilizes; on the other hand, if $v(f(s_n))$ is eventually strictly increasing for some $f\in K[X]$, we say that $E$ is of \emph{algebraic type}.
\end{Def}
The main difference between these two kinds  of sequences is the nature of the pseudo-limits: if $E$ is of algebraic type, then it has pseudo-limits in the algebraic closure $\overline{K}$ (for some extension $u$ of $v$), while if $E$ is of transcendental type then it admits a pseudo-limit only in a transcendental extension \cite[Theorems 2 and 3]{Kap}.

The central point of \cite{PS1} is the following: if $E=\{s_n\}_{n\in\N}\subset K$ is a pseudo-convergent sequence, then the set
\begin{equation}\label{V_E}
V_{E}=\{\phi\in K(X) \mid \phi(s_n)\in V,\textnormal{ for all but finitely many }n\in\N\}
\end{equation}
is a valuation domain of $K(X)$ extending $V$ \cite[Theorem 3.8]{PS1}. If $E,F$ are pseudo-convergent sequences of algebraic type, then $V_E=V_F$ if and only if $\limiti_E^u=\limiti_F^u$ for some extension $u$ of $v$ to $\overline{K}$ \cite[Theorem 5.4]{PS1}. In general, we say that two pseudo-convergent sequences $E,F$ are \emph{equivalent} if $V_E=V_F$; this condition can also be expressed by means of a notion analogous to the one defined classically for Cauchy sequences (see \cite[Definition 5.1]{PS1}).

We are interested in the study of the following subspace of $\Zar(K(X)|V)$:
$$\insiemeVE=\{V_E\mid E\subset K \textnormal{ is a pseudo-convergent sequence}\}.$$
The space $\insiemeVE$ is always regular under both the Zariski and the constructible topologies \cite[Theorem 6.15]{PS1}; however, these two topologies coincide if and only if the residue field of $V$ is finite \cite[Proposition 6.11]{PS1}.

\section{Fixed breadth}\label{sect:radfix}
In this section, we study the subsets of $\insiemeVE$ obtained by fixing the breadth of the pseudo-convergent sequences.
\begin{Def}\label{VEdelta}
Let $\delta\inR\cup\{+\infty\}$. We  denote by $\insiemeVE(\bullet,\delta)$ the set of valuation domains $V_E$ such that the breadth of $E$ is $\delta$.
\end{Def}

If $\delta=\infty$, then the elements of $\insiemeVE(\bullet,\delta)$ are the rings defined through pseudo-convergent sequences with $\Br(E)=(0)$, i.e., from pseudo-convergent sequences that are also Cauchy sequences. In this case, $E$ has a unique limit $\alpha\in\widehat{K}$, and by \cite[Remark 3.10]{PS1} we have
\begin{equation*}
V_E=W_\alpha=\{\phi\in K(X)\mid v(\phi(\alpha))\geq0\}.
\end{equation*}
Therefore, there is a natural bijection between $\widehat{K}$ and $\insiemeVE(\bullet,\infty)$, given by $\alpha\mapsto W_\alpha$; by \cite[Theorem 3.4]{PerTransc}, such a bijection is also a homeomoprhism, when $\widehat{K}$ is endowed with the $v$-adic topology and $\insiemeVE(\bullet,\infty)$ with the Zariski topology. In particular, it follows that the latter is an ultrametric space. Note that when $V$ is a discrete valuation ring, $\insiemeVE=\insiemeVE(\bullet,\infty)$.
\begin{prop}\label{prop:insVE-DVR}
Let $V$ be a discrete valuation ring. Then, $\insiemeVE\simeq\widehat{K}$ is an ultrametric space.
\end{prop}
\begin{proof}
The claim follows from the previous discussion and the fact that if $V$ is discrete then every pseudo-convergent sequence has infinite breadth.
\end{proof}

The purpose of this section is to see how the homeomorphism $\insiemeVE(\bullet,\infty)\simeq\widehat{K}$ generalizes when we consider pseudo-convergent sequences with fixed breadth $\delta\in\R$.

Fix $\delta\inR\cup\{\infty\}$, and set $r=e^{-\delta}$. Given two pseudo-convergent sequences $E=\{s_n\}_{n\inN}$ and $F=\{t_n\}_{n\inN}$, with $V_E,V_F\in\insiemeVE(\bullet,\delta)$, we set
\begin{equation*}
d_\sq(V_E,V_F)=\lim_{n\to\infty}\max\{d(s_n,t_n)-r,0\}.
\end{equation*}
It is clear that if $r=0$ (or, equivalently, $\delta=+\infty$) then $d_\sq(V_E,V_F)=d(\alpha,\beta)$, where $\alpha$ and $\beta$ are the (unique) limits of $E$ and $F$, respectively; so in this case we get the same distance as in \cite{PerTransc}. We shall interpret $d_\sq$ in a similar way in Proposition \ref{prop:dseqVEalgebraiccase}; we first show that it is actually a distance.
\begin{prop}\label{prop:dsq}
Preserve the notation above.
\begin{enumerate}[(a)]
\item\label{prop:dsq:wd} $d_\sq$ is well-defined.
\item\label{prop:dsq:dist} $d_\sq$ is an ultrametric distance on $\insiemeVE(\bullet,\delta)$.
\end{enumerate}
\end{prop}
\begin{proof}
\ref{prop:dsq:wd} Let $E=\{s_n\}_{n\inN}$ and $F=\{t_n\}_{n\inN}$ be two pseudo-convergent sequences. We start by showing that the limit of $a_n=\max\{d(s_n,t_n)-r,0\}$ exists. If all subsequences of $\{a_n\}_{n\inN}$ go to zero, we are done. Otherwise, there is a subsequence $\{a_{n_k}\}_{k\inN}$ with a positive (possibly infinite) limit; in particular, there is a $\overline{\delta}<\delta$ and $k_0\in\N$ such that $v(s_{n_k}-t_{n_k})<\overline{\delta}$ for all $k\geq k_0$. Choose $k_1\in\N$ such that $\overline{\delta}<\min\{\delta_{k_1},\delta'_{k_1}\}$ (where $\{\delta_n\}_{n\inN}$ and $\{\delta'_n\}_{n\inN}$ are the gauges of $E$ and $F$, respectively). Fix an $m=n_l$ such that $m>k_1$ and $l>k_0$. Then, for all $n> m$, we have
\begin{equation*}
v(s_n-t_n)=v(s_n-s_m+s_m-t_m+t_m-t_n)=v(s_m-t_m)
\end{equation*}
since $v(s_n-s_m)=\delta_m>\delta_{k_1}>\overline{\delta}>v(s_{n_l}-t_{n_l})=v(s_m-t_m)$, and likewise for $v(t_n-t_m)$. Hence, $a_n$ is eventually constant (more precisely, equal to $e^{-v(s_m-t_m)}-e^{-\delta}$); in particular, $\{a_n\}_{n\inN}$ has a limit.

In order to show that $d_{\delta}$ is well-defined, we need to show that, if $V_E=V_{E'}$, where $E=\{s_n\}_{n\inN}$ and $E'=\{s'_n\}_{n\inN}$, then 
\begin{equation*}
\lim_{n\to\infty}\max\{d(s_n,t_n)-r,0\}=\lim_{n\to\infty}\max\{d(s'_n,t_n)-r,0\}.
\end{equation*}
Let $l$ be the limit on the left hand side and $l'$ the limit on the right hand side. 

If $F$ is equivalent to $E$ and $E'$, by \cite[Definition 5.1 and Theorem 5.4]{PS1} for every $k$ there are $i_0,j_0,i'_0,j'_0$ such that $v(s_i-t_j)>\delta_k$, $v(s'_{i'}-t'_{j'})>\delta'_k$ for $i\geq i_0$, $j\geq j_0$, $i'\geq i'_0$, $j'\geq j'_0$. Hence, both $l$ and $l'$ are equal to 0, and in particular they are equal.

Suppose that $F$ is not equivalent to $E$ and $E'$. If $l$ is positive, and $\eta=-\log(l+r)$, then $v(s_n-t_n)=\eta$ for large $n$, and $\eta<\delta_k$ for some $k$; since $E$ and $E'$ are equivalent there is a $i_0$ such that $v(s_i-s'_i)>\delta_k$ for all $i\geq i_0$. Hence, for all large $n$,
\begin{equation*}
v(s'_n-t_n)=v(s'_n-s_n+s_n-t_n)=v(s_n-t_n)=\eta,
\end{equation*}
as claimed. The same reasoning applies if $l'>0$; furthermore, if $l=0=l'$ then clearly $l=l'$. Hence, $l=l'$ always, as claimed.

\ref{prop:dsq:dist} $d_\sq$ is obviously symmetric. Clearly $d_\sq(V_E,V_E)=0$; if $d_\sq(V_E,V_F)=0$, for every $r_k=e^{-\delta'_k}<r$ (where $\delta'_k=v(t_{k+1}-t_k)$) there is $i_0$ such that $d(s_i,t_i)<r_k$ for all $i\geq i_0$. Thus, if $i,j\geq i_0$, then
\begin{equation*}
d(s_i,t_j)=\max\{d(s_i,t_i),d(t_i,t_j)\}=r_k.
\end{equation*}
Hence, $E$ and $F$ are equivalent and $V_E=V_F$. The strong triangle inequality follows from the fact that $d(s_n,t_n)\leq\max\{d(s_n,s'_n),d(s'_n,t_n)\}$ for all $s_n,s'_n,t_n\in K$. Therefore, $d_\sq$ is an ultrametric distance.
\end{proof}

Let $\insiemeVE_K(\bullet,\delta)$ be the subset of $\insiemeVE(\bullet,\delta)$ corresponding to pseudo-convergent sequences with a pseudo-limit in $K$. We recall that by \cite[Theorem 5.4]{PS1} the map $V_E\mapsto \mathcal{L}_E$, from  $\insiemeVE_K(\bullet,\delta)$ to the set of closed balls in $K$ of radius $e^{-\delta}$, is a one-to-one correspondence. When $\delta=\infty$, $\insiemeVE_K(\bullet,\infty)$ corresponds to $K$ under the homeomorphism between $\insiemeVE(\bullet,\infty)$ and $\widehat{K}$; in particular, $\insiemeVE(\bullet,\infty)$ is the completion of $\insiemeVE_K(\bullet,\infty)$ under $d_\infty$. An analogous result holds for $\delta\inR$.

\begin{prop}\label{prop:compl} Let $\delta\in\R$. Then  
$\insiemeVE(\bullet,\delta)$ is the completion of $\insiemeVE_K(\bullet,\delta)$ under the metric $d_\sq$. In particular, $\insiemeVE(\bullet,\delta)$, under $d_\sq$, is a complete metric space.
\end{prop}
\begin{proof}
Let $\{\zeta_k\}_{k\inN}\subset\Gamma$ be an increasing sequence of real numbers with limit $\delta$ and, for every $k$, let $z_k$ be an element of $K$ of valuation $\zeta_k$; let $Z=\{z_k\}_{k\inN}$. It is clear that $Z$ is a pseudo-convergent sequence with $0$ as  a pseudo-limit and having breadth $\delta$. Then, for every $s\in K$, $s+Z=\{s+z_k\}_{k\inN}$ is a pseudo-convergent sequence with pseudo-limit $s$ and breadth $\delta$.

Let $E=\{s_n\}_{n\inN}$ be a pseudo-convergent sequence  with breadth $\delta$, and let $F_n=s_n+Z$. By above, $V_{F_n}\in \insiemeVE_K(\bullet,\delta)$, for each $n\in \N$. We claim that $\{V_{F_n}\}_{n\in\N}$ converges to $V_E$ in $\insiemeVE(\bullet,\delta)$. Indeed, fix $t\in\N$, and take $k>t$ such that $\zeta_k>\delta_t$. Then,
\begin{equation*}
u(s_t+z_k-s_k)=u(s_t-s_k+z_k)=\delta_t;
\end{equation*}
hence, $d(V_E,V_{F_n})=e^{-\delta_n}-e^{-\delta}$. In particular, the distance goes to 0 as $n\to\infty$, and thus $V_E$ is the limit of $V_{F_n}$.

Conversely, let $\{V_{F_n}\}_{n\inN}$ be a Cauchy sequence in $\insiemeVE_K(\bullet,\delta)$, and let $s_n\in K$ be a pseudo-limit of $F_n$. Then, $s_n+Z$ is another pseudo-convergent sequence with limit $s_n$ and breadth $\delta$; by \cite[Theorem 5.4]{PS1} it follows that $V_{F_n}=V_{s_n+Z}$. There is a subsequence of $E=\{s_n\}_{n\inN}$ which is pseudo-convergent; indeed, it is enough to take $\{s_{n_k}\}_{k\inN}$ such that $d(s_{n_k},s_{n_{k+1}})<d(s_{n_{k-1}},s_{n_k})$. Hence, without loss of generality $E$ itself is pseudo-convergent; we claim that $V_E$ is a limit of $\{V_{F_n}\}_{n\inN}$. Indeed, as above, $u(s_t+z_k-s_k)=\delta_t$ for large $k$, and thus $d_\sq(V_E,V_{s_n+Z})=e^{-\delta_t}-e^{-\delta}$. Thus, $\{V_{F_n}\}_{n\inN}$ has a limit, namely $V_E$. Therefore, $\insiemeVE(\bullet,\delta)$ is the completion of $\insiemeVE_K(\bullet,\delta)$.
\end{proof}

We now prove that the topology induced by $d_\sq$ is actually the Zariski topology.

\begin{teor}\label{teor:Sigmadelta}
Let $\delta\in\R\cup\{\infty\}$. On $\insiemeVE(\bullet,\delta)$, the Zariski topology, the constructible topology and the topology induced by $d_\sq$ coincide.
\end{teor}
\begin{proof}
If $\delta=\infty$, then the Zariski topology and the topology induced by $d_{\sq}$ coincide by \cite[Theorem 3.4]{PerTransc}.

Suppose now that $V$ is nondiscrete and fix $\delta\inR$. Let $V_E\in \insiemeVE(\bullet,\delta)$ and $\rho\in\R$, $\rho>0$: we show that the open ball $\ball(V_E,\rho)=\{V_F\in \insiemeVE(\bullet,\delta)\mid d_\sq(V_E,V_F)<\rho\}$  of the ultrametric topology induced by $d_\sq$ is open in the Zariski topology. Since  by Proposition \ref{prop:compl} $\insiemeVE_K(\bullet,\delta)$ is dense in $\insiemeVE(\bullet,\delta)$ under the metric $d_\sq$, without loss of generality we may assume that $V_E\in\insiemeVE_K(\bullet,\delta)$, i.e.,  $E$ has a pseudo-limit $b$ in $K$. To ease the notation, we denote by $B(\phi)$ the intersection $B(\phi)\cap\insiemeVE(\bullet,\delta)$.

Let $\gamma<\delta$ be such that $\rho=e^{-\gamma}-e^{-\delta}$. We claim that 
\begin{equation*}
\ball(V_E,\rho)=\bigcup_{\delta> v(c)>\gamma}B\left(\frac{X-b}{c}\right).
\end{equation*}
Indeed, suppose $V_F\in\ball(V_E,\rho)$, where $F=\{t_n\}_{n\inN}$. If $F$ is equivalent to $E$ then $V_E=V_F$ and $v\left(\frac{t_n-b}{c}\right)=\delta_n-v(c)$; since $\gamma<\delta$ and $\Gamma$ is dense in $\insR$, there is a $c\in K$ such that $\gamma<v(c)<\delta$, and for such a $c$ the limit of $\delta_n-v(c)$ is positive; hence, $V_E$ belongs to the union.  If $F$ is not equivalent to $E$, then $0<d_\sq(V_E,V_F)<\rho$, that is, $e^{-\delta}<\lim_nd(s_n,t_n)<e^{-\delta}+\rho$. By the proof of Proposition \ref{prop:dsq}\ref{prop:dsq:wd}, $v(s_n-t_n)$ is eventually constant, and thus there is an $\epsilon>0$ such that $\delta>v(s_n-t_n)\geq\gamma+\epsilon$ for all large $n$. Let $c\in K$ be of value comprised between $\gamma$ and $\gamma+\epsilon$ (such a $c$ exists because $\Gamma$ is dense in $\insR$); then,
\begin{equation*}
v\left(\frac{t_n-b}{c}\right)=v(t_n-b)-v(c)=v(t_n-s_n+s_n-b)-v(c)\geq\min\{\gamma+\epsilon,\delta_n\}-v(c)>0
\end{equation*}
since $\delta_n$ becomes bigger than $\gamma+\epsilon$. Hence, $\frac{X-b}{c}\in V_F$, or equivalently $V_F\in B\left(\frac{X-b}{c}\right)$. 

Conversely, suppose $V_F\neq V_E$ belongs to $B\left(\frac{X-b}{c}\right)$ for some $c\in K$ such that $\gamma<v(c)<\delta$. Since $\limiti_E\cap\limiti_F=\emptyset$ by \cite[Theorem 5.4]{PS1}, $b$ is not a pseudo-limit of $F$; therefore, $v(t_n-s_n)=v(t_n-b+b-s_n)=v(b-t_n)\geq v(c)>\gamma$ for sufficiently large $n$. Thus,
\begin{equation*}
d_\sq(V_E,V_F)=\lim_n d(s_n,t_n)-e^{-\delta}=\lim_n d(b,t_n)-e^{-\delta}<e^{-\gamma}-e^{-\delta}=\rho,
\end{equation*}
i.e., $V_F\in\ball(V_E,\rho)$. Thus, being the union of sets that are open in the Zariski topology, $\ball(V_E,\rho)$ is itself open in the Zariski topology. Therefore, the ultrametric topology is finer than the Zariski topology.

Let now $\delta$ be arbitrary, $\phi\in K(X)$ be a rational function, and suppose $V_E\in B(\phi)$ for some $V_E\in\insiemeVE(\bullet,\delta)$. We want to show that for some $\rho>0$ there is a ball $\ball(V_E,\rho)\subseteq B(\phi)$,  and thus that $B(\phi)$ is open in the ultrametric topology induced by $d_\sq$. We distinguish two cases.

Suppose that $E$ is of algebraic type, and let $\beta\in\limiti_E^u$ for some extension $u$ of $v$ to $\overline{K}$. By \cite[Lemma 6.6]{PS1}, there is an annulus $C=\corona(\beta,\tau,\delta)=\{s\in \overline{K}\mid \tau<u(s-\beta)<\delta\}$ such that $\phi(s)\in V$ for every $s\in C$.  Let $\epsilon=e^{-\tau}-e^{-\delta}$. Let $F=\{t_n\}_{n\inN}$ be a pseudo-convergent sequence with $d_\sq(V_E,V_F)<\epsilon$. Then, for every $n$ such that $e^{-\delta_n}-e^{-\delta}>d_\sq(V_E,V_F)$, we have
\begin{equation*}
d(t_n,\beta)=\max\{d(t_n,s_n),d(s_n,\beta)\}=e^{-\delta_n},
\end{equation*}
and in particular $v(t_n-\beta)$ becomes larger than $\tau$. Hence, $t_n$ is eventually in $C$ and $\phi(t_n)\in V$ for all large $n$, and thus $\phi\in V_F$; therefore, $\ball(V_E,\epsilon)\subseteq B(\phi)$.

Suppose that $E$ is of transcendental type. Let $\phi(X)=c\prod_{i=1}^A(X-\alpha_i)^{\epsilon_i}$ over $\overline K$, where each $\epsilon_i$ is either $1$ or $-1$. Then, there is an $N$ such that $u(s_n-\alpha_i)$ is constant for every $i$ and every $n\geq N$. Let $\delta'$ be the maximum among such constants; then, $\delta'<\delta$ (otherwise the $\alpha_i$ where such maximum is attained would be a pseudo-limit of $E$, in contrast to the fact that $E$ is of transcendental type). Let $\epsilon$ be such that $e^{-\delta}+\epsilon<e^{-\delta'}$ and let $V_F\in\ball(V_E,\epsilon)$, with $F=\{t_n\}_{n\inN}$. For all $i$, and all large $n$,
\begin{equation*}
d(t_n,\alpha_i)=\max\{d(t_n,s_n),d(s_n,\alpha_i)\}=d(s_n,\alpha_i),
\end{equation*}
and thus $u(t_n-\alpha_i)=u(s_n-\alpha_i)$. It follows that $v(\phi(t_n))=v(\phi(s_n))$ for large $n$; in particular, $v(\phi(t_n))$ is positive, and $\phi\in V_F$. Hence, $\ball(V_E,\epsilon)\subseteq B(\phi)$.

Hence, $B(\phi)$ is open under the topology induced by $d_\sq$ and therefore the Zariski topology and the topology induced by $d_\sq$ on $\insiemeVE(\bullet,\delta)$ are the same.

In order to prove that these topologies coincide also with the constructible topology, we need only to show that every $B(\phi)$, $\phi\in K(X)$, is closed in the Zariski topology. Let then $V_E\notin B(\phi)$. If $E$ is of transcendental type, exactly as above there exists $\epsilon>0$ such that for each $V_F\in \ball(V_E,\epsilon)$, where $F=\{t_n\}_{n\inN}$,  $v(\phi(t_n))=v(\phi(s_n))$ for large $n$; in particular, $v(\phi(t_n))$ is negative, and $\phi\notin V_F$; thus $\ball(V_E,\epsilon)$ is disjoint from $B(\phi)$. If $E$ is of algebraic type, then by \cite[Remark 6.7]{PS1}, there exists an annulus $\corona=\corona(\beta,\tau,\delta)$ such that $\phi(s)\notin V$ for every $s\in \corona$. As above, for every pseudo-convergent sequence $F=\{t_n\}_{n\in\N}$ with $d_\sq(V_E,V_F)<\epsilon$, with $\epsilon=e^{-\tau}-e^{-\delta}$, we have $t_n\in \corona$ for all but finitely many $n\in\N$, so that $\phi(t_n)\notin V$. Again, this shows that $\ball(V_E,\epsilon)$ is disjoint from $B(\phi)$.
\end{proof}

Joining Proposition \ref{prop:compl} with Theorem \ref{teor:Sigmadelta}, we obtain that the set $\insiemeVE_K=\{V_E\in\insiemeVE\mid \limiti_E\cap K\neq\emptyset\}=\bigcup_\delta\insiemeVE_K(\bullet,\delta)$ of all the extensions arising from pseudo-convergent sequences with pseudo-limits in $K$ is dense in $\insiemeVE$, with respect to both the Zariski and the constructible topology. This result can also be obtained as a corollary of \cite[Proposition 6.9]{PS1}.

If we restrict to pseudo-convergent sequences of algebraic type, the distance $d_\sq$ can be interpreted in a different way.
\begin{prop}\label{prop:dseqVEalgebraiccase}
Let $E,F\subset K$ be pseudo-convergent sequences of algebraic type with breadth $\delta$, and let $u$ be an extension of $v$ to $\overline{K}$. If $\beta\in\limiti_E^u$ and $\beta'\in\limiti_F^u$, then 
\begin{equation*}
d_\sq(V_E,V_F)=\max\{d_u(\beta,\beta')-e^{-\delta},0\}.
\end{equation*}
\end{prop}
\begin{proof}
If $d_u(\beta,\beta')\leq e^{-\delta}$, then the pseudo-limits of $E$ and $F$ coincide, and thus $V_E=V_F$ by \cite[Theorem 5.4]{PS1}; hence, $d_\sq(V_E,V_F)=0$. On the other hand, if $d_u(\beta,\beta')>e^{-\delta}$ then $u(\beta-\beta')<\delta$ and thus, for large $n$,
\begin{equation*}
v(s_n-t_n)=u(s_n-\beta+\beta-\beta'+\beta'-t_n)=u(\beta-\beta');
\end{equation*}
hence, $d_\sq(V_E,V_F)=d_u(\beta,\beta')-e^{-\delta}$, as claimed. 
\end{proof}

If $V$ is a DVR, then $\mathcal{V}=\mathcal{V}(\bullet,\infty)$, so, in this case, the distance $d_\infty$ is an ultrametric distance on the whole $\mathcal{V}$. On the other hand, if $V$ is not discrete, it is not possible to unify the metrics $d_\sq$ in a single metric defined on the whole $\insiemeVE$. To this end, we need the following  lemma.
\begin{lemma}\label{lemma:chiusura-bdelta}
Let $\delta\in\R\cup\{\infty\}$. Then the closure of $\insiemeVE(\bullet,\delta)$ in $\insiemeVE$ is equal to $\displaystyle{\bigcup_{\delta'\leq\delta}\insiemeVE(\bullet,\delta')}$.
\end{lemma}
\begin{proof}
If $V$ is discrete, then the statement is a tautology (see Proposition \ref{prop:insVE-DVR}). We assume henceforth that $V$ is not discrete. 

Let $E=\{s_n\}_{n\inN}$ be a pseudo-convergent sequence with breadth $\delta'<\delta$; we want to show that $V_E$ is in the closure of $\insiemeVE(\bullet,\delta)$. By Proposition \ref{prop:compl}, $\insiemeVE(\bullet,\delta')$ is contained in the closure of $\insiemeVE_K(\bullet,\delta')$; hence, we can suppose that $E$ has a pseudo-limit in $K$.

For each $n\in\N$, let $E_n$ be a pseudo-convergent sequence with pseudo-limit $s_n$ and breadth $\delta$: since $\delta'<\delta$, by \cite[Proposition 6.9]{PS1} $V_E$ is the limit of $V_{E_n}$ in the Zariski topology, and thus it belongs to the closure of $\insiemeVE_K(\bullet,\delta')$, as claimed. If $\delta=\infty$ we are done; suppose for the rest of the proof that $\delta<\infty$.

Suppose $\delta'>\delta$; we claim that if $E=\{s_n\}_{n\inN}$ is pseudo-convergent sequence with breadth $\delta'$ then there is an open set containing $V_E$ and disjoint from $\insiemeVE(\bullet,\delta)$. Let $\gamma\in\Gamma_v$ be such that $\delta'>\gamma>\delta$; then, there is an $N$ such that $v(s_n-s_{n+1})>\gamma$ for all $n\geq N$. Take $s=s_N$, and consider the open set $B\left(\frac{X-s}{c}\right)$, where $c\in K$ has value $\gamma$. Then, $V_E\in B\left(\frac{X-s}{c}\right)$ since $v(s_n-s_N)=\delta'_N>\gamma$ for all $n\geq N$. On the other hand, if $F=\{t_n\}_{n\inN}\subset K$ is a pseudo-convergent sequence of breadth $\delta$  and $V_F\in B\left(\frac{X-s}{c}\right)$, then $F$ would be eventually contained in the ball of center $s$ and radius $\gamma$, and in particular $v(t_n-t_{n+1})\geq\gamma$ for all large $n$. However, $v(t_n-t_{n+1})<\delta<\gamma$, a contradiction. Therefore, $V_F\notin B\left(\frac{X-s}{c}\right)$ and so $V_E$ is not in the closure of $\insiemeVE(\bullet,\delta)$.
\end{proof}

\begin{prop}\label{prop:restriz-dsq}
Let $V$ be a rank one non-discrete valuation domain. Suppose $\insiemeVE$ is metrizable with a metric $d$. Then, for any $\delta\in\insR\cup\{\infty\}$, the restriction of $d$ to $\insiemeVE(\bullet,\delta)$ is \emph{not} equal to $d_\sq$.
\end{prop}
\begin{proof}
If the restriction of $d$ is equal to $d_\sq$, then by Proposition \ref{prop:compl} $\insiemeVE(\bullet,\delta)$ would be complete with respect to $d$. However, this would imply that $\insiemeVE(\bullet,\delta)$ is closed, in contrast to Lemma \ref{lemma:chiusura-bdelta}.
\end{proof}

To conclude this section, we analyze the relationship among  the sets $\insiemeVE(\bullet,\delta)$, as $\delta$ ranges in $\R\cup\{\infty\}$. Recall that two metric spaces $(X,d)$ and $(X',d')$ are \emph{similar} if there is a map $\psi:X\longrightarrow X'$ and a constant $r>0$ such that $d'(\psi(x),\psi(y))=rd(x,y)$ for every $x,y\in X$. We call such a map $\psi$ a \emph{similitude}.
\begin{prop}
If $\delta_1-\delta_2\in\Gamma_v$, then the metric spaces $(\insiemeVE(\bullet,\delta_1),d_{\delta_1})$ and $(\insiemeVE(\bullet,\delta_2),d_{\delta_2})$ are similar; in particular, they are homeomorphic when endowed with the Zariski topology.
\end{prop}
\begin{proof}
Given a pseudo-convergent sequence $E=\{s_n\}_{n\inN}$ and $c\in K$, $c\neq 0$, we denote by $cE$ the sequence $\{cs_n\}_{n\inN}$. Clearly, $cE$ is again pseudo-convergent, it has breadth $\delta_E+v(c)$, and two sequences $E$ and $F$ are equivalent if and only if $cE$ and $cF$ are equivalent.

Let $c\in K$ be such that $v(c)=\delta_1-\delta_2$. Then, the map
\begin{equation*}
\begin{aligned}
\Psi_c\colon \insiemeVE(\bullet,\delta_2) & \longrightarrow \insiemeVE(\bullet,\delta_1)\\
V_E & \longmapsto V_{cE}
\end{aligned}
\end{equation*}
is well-defined and bijective (its inverse is $\Psi_{c^{-1}}:\insiemeVE(\bullet,\delta_1)\longrightarrow\insiemeVE(\bullet,\delta_2)$). We claim that $\Psi_c$ is a similitude. Indeed, let $E=\{s_n\}_{n\inN}$ and $F=\{t_n\}_{n\inN}$ be pseudo-convergent sequences of breadth $\delta_2$, and suppose $V_E\neq V_F$. By the proof of Proposition \ref{prop:dsq}, there is an $N$ such that $v(s_n-t_n)=v(s_N-t_N)$ for all $n\geq N$. Hence, for these $n$'s,
\begin{align*}
e^{-v(cs_n-ct_n)}-e^{-\delta_1}=e^{-v(c)}e^{-v(s_n-t_n)}-e^{-\delta_1}=e^{-v(c)}[e^{-v(s_n-t_n)}-e^{-\delta_2}]
\end{align*}
so that, passing to the limit,  $d_{\delta_1}(V_{cE},V_{cF})=e^{-v(c)}d_{\delta_2}(V_E,V_F)$. Hence, $\Psi_c$ is an similitude, and in particular a homeomorphism when $\insiemeVE(\bullet,\delta_1)$ and $\insiemeVE(\bullet,\delta_2)$ are endowed with the metric topology. Since this topology coincides with the Zariski topology (Theorem \ref{teor:Sigmadelta}), they are homeomorphic also under the Zariski topology.
\end{proof}

\section{Fixed pseudo-limit}\label{sect:bfix}
In the previous section, we considered valuation domains induced by pseudo-convergent sequences having the same breadth; in this section, we reverse the situation by considering pseudo-convergent sequences having a prescribed pseudo-limit. Note that, in particular, these pseudo-convergent sequences are of algebraic type.

Throughout this section, let $u$ be a fixed extension of $v$ to $\overline K$. 
\begin{Def}
Let $\beta\in \overline K$. We set
$$\insiemeVE^u(\beta,\bullet)=\{V_E\in\insiemeVE \mid \beta\in\limiti_E^u\}$$
To ease the notation, we set $\insiemeVE^u(\beta,\bullet)=\insiemeVE(\beta,\bullet)$.
\end{Def}
Equivalently, a valuation domain $V_E$ is in $\insiemeVE(\beta,\bullet)$ if $\beta$ is a center of $\limiti_E^u$, i.e., if $\limiti_E^u=\{x\in \overline K \mid u(x-\beta)\geq\delta_E\}$. Note that if $V_E\in\insiemeVE^u(\beta,\bullet)$ then $E$ must be of algebraic type, since it must have a pseudo-limit in $\overline{K}$.

If $V$ is a DVR, then $\insiemeVE(\beta,\bullet)$ reduces to the single element $W_{\beta}=\{\phi\in K(X)\mid \phi(\beta)\in V\}$ (see  \cite[Remark 3.10]{PS1}), which corresponds to any Cauchy sequence $E\subset K$  converging to $\beta$.

We start by showing that each $\insiemeVE(\beta,\bullet)$ is closed in $\insiemeVE$.
\begin{prop}\label{prop:betachiuso}
Let $\beta\in\overline{K}$, and let $u$ be an extension of $v$ to $\overline{K}$. Then, $\insiemeVE(\beta,\bullet)=\insiemeVE^u(\beta,\bullet)$ is closed in $\insiemeVE$ with respect to the Zariski topology.
\end{prop}
\begin{proof}
If $V$ is discrete, then $\insiemeVE(\beta,\bullet)$ has just one element (see the comments above). By \cite[Theorem 3.4]{PerTransc} each point of $\insiemeVE$ is closed, so the statement is true in this case. Henceforth, for the rest of the proof we assume that $V$ is non discrete.

Let $V_E\notin\insiemeVE(\beta,\bullet)$. We distinguish two cases.

Suppose first that $E=\{s_n\}_{n\inN}$ is of algebraic type, and let $\alpha\in\overline{K}$ be a pseudo-limit of $E$ with respect to $u$. Since $\beta\notin \mathcal L_E\Leftrightarrow u(\alpha-\beta)<\delta_E$ (Lemma \ref{lemma:kaplansky}) it follows that there is $m\in\N$ such that $u(\alpha-\beta)<u(\alpha-s_m)$. Let $s=s_m$. Choose a $d\in K$ such that 
\begin{equation*}
u(\beta-\alpha)=u(\beta-s)<v(d)<u(\alpha-s)<\delta_E,
\end{equation*}
and let $\phi(X)=\frac{X-s}{d}$; we claim that $V_E\in B(\phi)$ but $B(\phi)\cap\insiemeVE(\beta,\bullet)=\emptyset$.

Indeed,
\begin{equation*}
v(\phi(s_n))=v\left(\frac{s_n-s}{d}\right)=v(s_n-s)-v(d)>0
\end{equation*}
since $v(s_n-s)=u(s_n-\alpha+\alpha-s)=u(\alpha-s)$ for large $n$; hence $V_E\in B(\phi)$. On the other hand, if $F=\{t_n\}_{n\inN}$ has pseudo-limit $\beta$, then $v(t_n-s)=u(t_n-\beta+\beta-s)=u(\beta-s)$ for large $n$ and so
\begin{equation*}
v(\phi(t_n))=u(\beta-s)-v(d)<0,
\end{equation*}
i.e., $V_F\notin B(\phi)$. The claim is proved.

Suppose now that $E=\{s_n\}_{n\inN}$ is of transcendental type: then, $u(s_n-\beta)$ is eventually constant, say equal to $\lambda$. Then, $\lambda<\delta$, for otherwise $\beta$ would be a pseudo-limit of $E$; hence, we can take a $d\in K$ such that $\lambda<v(d)<\delta$. Choose an $N$ such that $u(s_N-\beta)=\lambda$ and such that $v(d)<\delta_N$, and define $\phi(X)=\frac{X-s_N}{d}$. Then, $v(\phi(s_n))=\delta_N-v(d)>0$ for $n>N$, and thus $V_E\in B(\phi)$. Suppose now $v(\phi(t))\geq 0$. Then, $v(t-s_N)\geq v(d)>\lambda$; however, $v(t-s_N)=u(t-\beta+\beta-s_N)$, and since $u(\beta-s_N)=\lambda$ we must have $u(t-\beta)=\lambda$. In particular, there is no annulus $C$ of center $\beta$ such that $\phi(t)\in V$ for all $t\in C$; hence, by \cite[Lemma 6.6]{PS1}, $V_F\notin B(\phi)$ for every $V_F\in\insiemeVE(\beta,\bullet)$, i.e., $\insiemeVE(\beta,\bullet)\cap B(\phi)=\emptyset$. The claim is proved.
\end{proof}

We now want to characterize the Zariski topology of $\insiemeVE(\beta,\bullet)$.  By \cite[Theorem 5.4]{PS1}, there is a natural injective map
\begin{equation}\label{SigmaBeta}
\begin{aligned}
\Sigma_\beta\colon\insiemeVE(\beta,\bullet) & \longrightarrow(-\infty,+\infty] \\
V_E & \longmapsto \delta_E.
\end{aligned}
\end{equation}
In general this map is not surjective: for example, there might be some $\beta\in\overline K$ which is not the limit of any Cauchy sequence in $K$ (with respect to $u$) and thus $\delta_E\neq+\infty$ for every $V_E\in\insiemeVE(\beta,\bullet)$. By \cite[Proposition 5.5]{PS1} the image of $\Sigma_\beta$ is $(-\infty,\delta(\beta,K)]$, where $\delta(\beta,K)$ is defined as
\begin{equation*}
\delta(\beta,K)=\sup\{u(\beta-x)\mid x\in K\}.
\end{equation*}

In order to study the Zariski topology on $\mathcal V(\beta,
\bullet)$, we introduce a topology on the interval $(-\infty,\delta(\beta,K)]$.
\begin{Def}\label{def:lowerlimitGamma}
Let $a,b\inR\cup\{-\infty,+\infty\}$, with $a<b$, and let $\Lambda\subseteq\insR$. The \emph{$\Lambda$-upper limit topology} on $(a,b]$ is the topology generated by the sets $(\alpha,\lambda]$, for $\lambda\in\Lambda\cup\{\infty\}$ and $\alpha\in(a,b]$. We denote this space by $\limtop{a}{b}{\Lambda}$.
\end{Def}
The $\Lambda$-upper limit topology is a variant of the upper limit topology (see e.g. \cite[Counterexample 51]{counterexamples-topology}), and in fact the two topologies coincide when $\Lambda=\insR$.

For the next theorem we need to recall a definition and a result from \cite{PS1}. Let $E=\{s_n\}_{n\inN}$ be a pseudo-convergent sequence; we can associate to $E$ the map
\begin{equation*}
\begin{aligned}
w_E\colon K(X) & \longrightarrow \insR\cup\{\infty\}\\
\phi & \longmapsto \lim_{n\to\infty}v(\phi(s_n));
\end{aligned}
\end{equation*}
this map is always well-defined, and it is possible to characterize when it is a valuation on $K(X)$ \cite[Propositions 4.3 and 4.4]{PS1}. Given $s\in K$ and $\gamma\inR$, we set
\begin{equation*}
\Omega(s,\gamma)=\{V_F\in\insiemeVE \mid w_F(X-s)\leq\gamma\};
\end{equation*}
this set is always open and closed in $\insiemeVE$ (with respect to the Zariski topology) \cite[Lemma 6.14]{PS1}.

\begin{teor}\label{teor:betafix}
Suppose $V$ is not discrete, and let $\beta\in\overline{K}$ be a fixed element. The map $\Sigma_\beta$ defined in \eqref{SigmaBeta} is a homeomorphism between $\insiemeVE(\beta,\bullet)$ (endowed with the Zariski topology) and $\limtopstd$.
\end{teor}
\begin{proof}
To shorten the notation, let $\mathcal{X}=\limtopstd$.

We start by showing that $\Sigma_\beta$ is continuous. Clearly, $\Sigma_\beta^{-1}(\mathcal{X})=\insiemeVE(\beta,\bullet)$ is open. 

Suppose $\gamma\in\insQ\Gamma_v$ satisfies $\gamma<\delta(\beta,K)$. Then, there is a $t\in K$ such that $u(t-\beta)>\gamma$; we claim that
\begin{equation*}
\Sigma_\beta^{-1}((-\infty,\gamma])=\Omega(t,\gamma)\cap\insiemeVE(\beta,\bullet).
\end{equation*}
Indeed, let $E=\{s_n\}_{n\inN}$ be a pseudo-convergent sequence having $\beta$ as a pseudo-limit. If $\delta_E\leq\gamma$, then (since $u(\beta-t)>\gamma$)
\begin{equation*}
w_E(X-t)=\lim_{n\to\infty}v(s_n-t)=\lim_{n\to\infty}u(s_n-\beta+\beta-t)=\delta_E
\end{equation*}
and so $V_E\in\Omega(t,\gamma)$. Conversely, if $V_E\in\Omega(t,\gamma)\cap\insiemeVE(\beta,\bullet)$ then $w_E(X-t)\leq\gamma$, and thus (using again $u(\beta-t)>\gamma$)
\begin{equation*}
\delta_E=\lim_{n\to\infty}u(s_n-\beta)=\lim_{n\to\infty}u(s_n-t+t-\beta)=\lim_{n\to\infty}u(s_n-t)=w_E(X-t)\leq\gamma,
\end{equation*}
i.e., $\Sigma_\beta(V_E)\leq \gamma$.

By \cite[Lemma 6.14]{PS1}, $\Omega(t,\gamma)$ is open and closed in $\insiemeVE$; hence, $\Sigma_\beta^{-1}((-\infty,\gamma])$ and $\Sigma_\beta^{-1}((\gamma,\delta(\beta,K)])$ are both open. If now $(a,b]$ is an arbitrary basic open set of $\mathcal{X}$, with $b\in\insQ\Gamma$, then
\begin{equation*}
\Sigma_\beta^{-1}((a,b])=\Sigma_\beta^{-1}((-\infty,b])\cap\left(\bigcup_{\substack{c\in\Q\Gamma_v\\ c>a}}\Sigma_\beta^{-1}((c,\delta(\beta,K)])\right)
\end{equation*}
is open. Hence, $\Sigma_\beta$ is continuous.

Let now $\phi$ be an arbitrary nonzero rational function over $K$, and for ease of notation let $B(\phi)$ denote the intersection $B(\phi)\cap\insiemeVE(\beta,\bullet)$. Suppose $\delta\in\Sigma_\beta(B(\phi))$, and let $E=\{s_n\}_{n\inN}$ be a pseudo-convergent sequence of breadth $\delta$ having $\beta$ as a pseudo-limit. By \cite[Lemma 6.6]{PS1} there are $\theta_1,\theta_2\in\insQ\Gamma_v$ such that $\theta_2<\delta\leq\theta_1$ and such that $v(\phi(t))\geq 0$ for all $t\in\corona(\beta,\theta_1,\theta_2)$. In particular, if $V_F\in \insiemeVE(\beta,\bullet)$, $F=\{t_n\}_{n\in\N}$, is such that $\Sigma_{\beta}(V_F)\in(\theta_1,\theta_2]$ we have that $t_n\in\corona(\beta,\theta_1,\theta_2)$ for each $n\geq N$, for some $N\in\N$, so that $v(\phi(t_n))\geq0$ for each $n\geq N$, thus $\phi\in V_F$. Hence, $(\theta_1,\theta_2]\subseteq \Sigma_\beta(B(\phi))$, and thus $(\theta_1,\theta_2]$ is an open neighbourhood of $\delta$ in $\Sigma_\beta(B(\phi))$, which thus is open.

Hence, $\Sigma_\beta$ is open, and thus $\Sigma_\beta$ is a homeomorphism.
\end{proof}

Let $\mathcal W$ be the set of valuation domains of $K(X)$ associated to the valuations $w_E$ defined above, as $E$ ranges through the set of pseudo-convergent sequences of $K$ such that $w_E$ is a valuation. When $V$ is not discrete, we obtain a new proof of the non-compactness of $\mathcal{W}$, independent from \cite[Proposition 6.4]{PS1}.
\begin{cor}\label{cor:Vcompatto}
The spaces $\insiemeVE$ and $\mathcal{W}$ are not compact.
\end{cor}
\begin{proof}
If $V$ is a DVR, then $\mathcal{V}$ is homeomorphic to $\widehat{K}$ (\cite[Theorem 3.4]{PerTransc}). In particular, it is not compact. The space $\mathcal{W}$ is not compact by \cite[Proposition 6.4]{PS1}.

Suppose that $V$ is not discrete, and let $\beta\in\overline{K}$ be a fixed element. By Proposition \ref{prop:betachiuso}, $\insiemeVE(\beta,\bullet)$ is closed in $\insiemeVE$; hence if $\insiemeVE$ were compact so would be $\insiemeVE(\beta,\bullet)$. By Theorem \ref{teor:betafix}, it would follow that $\mathcal{X}=\limtopstd$ is compact. However, let $\gamma_1>\gamma_2>\cdots$ be a decreasing  sequence of elements in $\insQ\Gamma_v$, with $\delta(\beta,K)>\gamma_1$. Then, the family $(\gamma_1,\delta(\beta,K)],(\gamma_2,\gamma_1],\ldots,(\gamma_{n+1},\gamma_n],\ldots$ is an open cover of $\mathcal{X}$ without finite subcovers: hence, $\mathcal{X}$ is not compact, and so neither is $\insiemeVE$.

Let $\Psi:\mathcal{W}\longrightarrow\insiemeVE$ be the map $W_E\mapsto V_E$ (see \cite[Proposition 6.13]{PS1}). Since $\Psi$ is continuous, if $\mathcal{W}$ were compact then so would be its image $\insiemeVE_0$. Hence, as in the previous part of the proof, also $\insiemeVE_0\cap\insiemeVE(\beta,\bullet)$ would be compact; however, since $\Sigma_\beta(\insiemeVE_0\cap\insiemeVE(\beta,\bullet))=(-\infty,\delta(\beta,K)]\setminus\{+\infty\}$, we can use the same method as above (eventually substituting $(\gamma_1,+\infty]$ with $(\gamma_1,+\infty)$) to show that this set can't be compact. Hence, $\mathcal{W}$ is not compact, as claimed.
\end{proof}

We note that, when $V$ is a DVR, $\widehat{K}$ (and thus $\mathcal{V}$) is locally compact if and only if the residue field of $V$ is finite \cite[Chapt. VI, \S 5, 1., Proposition 2]{Bourb}. We conjecture that $\mathcal{V}$ is locally compact also when $V$ is not discrete.

\begin{prop}\label{prop:betafix-constr}
Let $\beta\in\overline{K}$, and let $u$ be an extension of $v$ to $\overline{K}$. Then, the Zariski and the constructible topologies agree on $\insiemeVE(\beta,\bullet)=\insiemeVE^u(\beta,\bullet)$.
\end{prop}
\begin{proof}
It is enough to show that $B(\phi)\cap\insiemeVE(\beta,\bullet)$ is closed for every $\phi\in K(X)$. Suppose $\delta\in C=\Sigma_\beta(\insiemeVE(\beta,\bullet)\setminus B(\phi))$ and let $V_E\in \insiemeVE(\beta,\bullet)\setminus B(\phi)$: by \cite[Lemma 6.6 and Remark 6.7]{PS1}, there is an annulus $\corona=\corona(\beta,\theta_1,\theta_2)$ with $\theta_1,\theta_2\in \insQ\Gamma_v$, $\theta_1<\delta\leq\theta_2$ and such that $\phi(t)\notin V$ for all $t\in\corona$.  Hence, $(\theta_1,\theta_2]$ is an open neighborhood of $\delta$ in $\limtopstd$ contained in $C$; thus, $C$ is open and $B(\phi)\cap\insiemeVE(\beta,\bullet)$ is closed, being the complement of the image of $C$ under the homeomorphism $\Sigma_{\beta}^{-1}$ (see Theorem \ref{teor:betafix}).
\end{proof}

To conclude, we study the metrizability of $\insiemeVE(\beta,\bullet)$ and $\insiemeVE$. It is well-known \cite[Counterexample 51(4)]{counterexamples-topology} that the upper limit topology is not metrizable, since it is separable but not second countable. Something similar happens for $\limtop{a}{b}{\Lambda}$.

\begin{prop}\label{prop:XGamma-metriz}
Let $\Lambda$ be a subset of $(a,b]$ that is dense in the Euclidean topology. The following are equivalent:
\begin{enumerate}[(i),noitemsep]
\item\label{prop:XGamma-metriz:num} $\Lambda$ is countable;
\item\label{prop:XGamma-metriz:2num} $\limtop{a}{b}{\Lambda}$ is second-countable;
\item\label{prop:XGamma-metriz:metriz} $\limtop{a}{b}{\Lambda}$ is metrizable;
\item\label{prop:XGamma-metriz:ultrametriz} $\limtop{a}{b}{\Lambda}$ is an ultrametric space.
\end{enumerate}
\end{prop}
\begin{proof}
\ref{prop:XGamma-metriz:metriz} $\Longrightarrow$\ref{prop:XGamma-metriz:2num} follows from the fact that $\limtop{a}{b}{\Lambda}$ is separable (since, for example, $\insQ\cap(a,b]$ is dense in $\limtop{a}{b}{\Lambda})$; \ref{prop:XGamma-metriz:ultrametriz} $\Longrightarrow$ \ref{prop:XGamma-metriz:metriz} is obvious.

\ref{prop:XGamma-metriz:2num} $\Longrightarrow$ \ref{prop:XGamma-metriz:num} Any basis of $\limtop{a}{b}{\Lambda}$ must contain an open set of the form $(\alpha,\lambda]$, for each $\lambda\in\Lambda$ (and some $\alpha\in(-\infty,\lambda)$). Hence, if $\limtop{a}{b}{\Lambda}$ is second-countable then $\Lambda$ must be countable.

\ref{prop:XGamma-metriz:num} $\Longrightarrow$ \ref{prop:XGamma-metriz:ultrametriz} Suppose that $\Lambda$ is countable, and fix an enumeration $\{\lambda_1,\lambda_2,\ldots\}$ of $\Lambda$. Let $r:\Lambda\longrightarrow\insR$ be the map sending $\lambda_i$ to $1/i$; then, for each $x,y\in(a,b]$ we set
\begin{equation*}
d(x,y)=\left\{
\begin{array}{ll}
\max\{r(\lambda)\mid \lambda\in[\min(x,y),\max(x,y))\cap\Lambda\},& \textnormal{ if }x\not=y\\
0,&\textnormal{ if }x=y.
\end{array}
\right.
\end{equation*}
We claim that $d$ is a metric on $(a,b]$ whose topology is exactly $\limtop{a}{b}{\Lambda}$.

Note first that $d$ is well-defined and nonnegative; it is also clear from the definition (and the fact that $\Lambda$ is dense in $\insR$) that $d(x,y)=0$ if and only if $x=y$ and that $d(x,y)=d(y,x)$. Let now $x,y,z\in(a,b]$, and suppose without loss of generality that $x\leq y$. If $z\leq x$, then $[z,y)\supseteq[x,y)$, and thus $d(x,y)\leq d(y,z)$; in the same way, if $y\leq z$ then $[x,z)\supseteq[x,y)$ and $d(x,y)\leq d(x,z)$. If $x\leq z\leq y$, then $[x,y)=[x,z)\cup[z,y)$; hence, $d(x,y)=\max\{d(x,z),d(y,z)\}$. In all cases, we have $d(x,y)\leq\max\{d(x,z),d(y,z)\}$, and thus $d$ induces an ultrametric space.

Let now $x\in\Lambda\subseteq(a,b]$ and $\rho\inR$ be positive; we claim that the open ball $B=B_d(x,\rho)=\{t\in(a,b]\mid d(x,t)<\rho\}$ is equal to $(y,z]$, where
\begin{align*}
y & =\max\{\lambda\in\Lambda\cap(-\infty,x)\mid r(\lambda)\geq \rho\},\\
z & =\min\{\lambda\in\Lambda\cap(x,+\infty)\mid r(\lambda)\geq \rho\}
\end{align*}
(with the convention $\max\emptyset=a$ and $\min\emptyset=b$). Note that since $\rho>0$, there are only a finite number of $\lambda$ with $r(\lambda)\geq\rho$; in particular, $y,z\in\Lambda$ and by definition, $y<x< z$.

Let $t\in(a,b]$. If $t<y$, then $r(\lambda)\geq \rho$ for some $\lambda\in(t,x)\cap\Lambda$, and thus $d(t,x)\geq\rho$, and so $t\notin B$; in the same way, if $y<t<x$, then $r(\lambda)<\rho$ for every $\lambda\in(t,x)\cap\Lambda$, and thus $t\in B$. Symmetrically, if $x<t<z$ then $t\in B$, while if $z<t$ then $t\notin B$. We thus need to analyze the cases $t=y$ and $t=z$.

By definition,
\begin{equation*}
d(x,z)=\max\{r(\lambda)\mid\lambda\in[x,z)\cap\Lambda\};
\end{equation*}
since by definition $r(\lambda)<\rho$ for every $\lambda\in[x,z)\cap\Lambda$, we have $d(x,z)<\rho$ and $z\in B_d(x,\rho)$.

Since $y\in\Lambda$, we have $r(y)\geq\rho$. Thus,
\begin{equation*}
d(x,y)=\max\{r(\lambda)\mid\lambda\in[y,x)\cap\Lambda\}\geq r(y)\geq \rho
\end{equation*}
and $y\notin B_d(x,\rho)$. Thus, $B_d(x,\rho)=(y,z]$ as claimed; therefore, $B_d(x,\rho)$ is open in $\limtop{a}{b}{\Lambda}$.

The family of the intervals $(y,z]$, as $z$ ranges in $\Lambda$ and $y$ in $(a,b]$, is a basis of $\limtop{a}{b}{\Lambda}$; therefore, the topology induced by $d$ on $(a,b]$ is exactly the $\Lambda$-upper limit topology. Hence, $\limtop{a}{b}{\Lambda}$ is an ultrametric space, as claimed.
\end{proof}

As a consequence, we obtain a necessary condition for metrizability, while in \cite[Corollary 6.16]{PS1} we obtained a sufficient condition, namely, if $V$ is countable, then $\mathcal V$ is metrizable.
\begin{cor}\label{cor:Gammav-uncont-metriz}
Let $V$ be a valuation ring with uncountable value group. Then, $\insiemeVE$ and $\Zar(K(X)|V)^\cons$ are not metrizable.
\end{cor}
\begin{proof}
If $\insiemeVE$ were metrizable, so would be $\insiemeVE(\beta,\bullet)$, in contrast to Theorem \ref{teor:betafix} and Proposition \ref{prop:XGamma-metriz} (note that, if the value group of $V$ is uncountable, in particular $V$ is not discrete). Similarly, if $\Zar(K(X)|V)^\cons$ were metrizable, so would be $\insiemeVE(\beta,\bullet)$, endowed with the constructible topology. Since the Zariski and the constructible topologies agree on $\insiemeVE(\beta,\bullet)$ (Proposition \ref{prop:betafix-constr}), this is again impossible.
\end{proof}

\section{Beyond pseudo-convergent sequences}\label{sect:beyond}
Corollary \ref{cor:Gammav-uncont-metriz} gives a condition for the non-metrizability of $\Zar(K(X)|V)^\cons$ that depends on the value group of $V$. In this section we prove a similar criterion, but based on the residue field of $V$.
\begin{lemma}\label{lemma:quoziente}
Let $V$ be a valuation ring with quotient field $K$, let $L$ be an extension field of $K$ and let $W$ be an extension of $V$ to $L$. Let $\pi:W\longrightarrow W/M$ be the quotient map. Then, the map
\begin{equation*}
\begin{aligned}
\{Z\in\Zar(L|V)\mid Z\subseteq W\} & \longrightarrow\Zar(W/M_W|V/M_V),\\
Z & \longmapsto \pi(Z)\\
\end{aligned}
\end{equation*}
is a homeomorphism, when both sets are endowed with either the Zariski or the constructible topology.
\end{lemma}
\begin{proof}
Apply \cite[Lemma 4.2]{isolated-Zariski} with $D=V$.
\end{proof}

\begin{lemma}\label{lemma:metriz-uncount-limit}
Let $X$ be an uncountable compact topological space with at most one limit point. Then, $X$ is not metrizable.
\end{lemma}
\begin{proof} Since $X$ is infinite and compact it has a limit point, say $x_0$, which is also unique by assumption.   
Suppose that $X$ is metrizable, and let $d$ be a metric inducing the topology. For each integer $n>0$, let $C_n=\{y\in X\mid 1/n\leq d(y,x_0)\}$. By construction, $x_0\notin C_n$, and thus all points of $C_n$ are isolated. Furthermore, $C_n$ is closed (since it is the complement of an open ball), and thus it is compact; therefore, $C_n$ must be finite. Hence, the countable union $\bigcup_{n>0}C_n$ is a countable set, against the fact that the union is equal to the uncountable set $X\setminus\{x_0\}$. Therefore, $X$ is not metrizable.
\end{proof}

\begin{prop}\label{prop:uncount-resfield}
Let $V$ be a valuation ring with uncountable residue field. Then, $\Zar(K(X)|V)^\cons$ is not metrizable.
\end{prop}
\begin{proof}
Let $W$ be the Gaussian extension of $V$ (see e.g. \cite{gilmer}); then, $W$ is an extension of $V$ to $K(X)$ having the same value group of $V$ and whose residue field is $k(t)$, where $k=V/M$ is the quotient field of $V$ and $t$ is an indeterminate. Consider $\Delta=\{Z\in\Zar(K(X)|V)\mid Z\subseteq W\}$; by Lemma \ref{lemma:quoziente}, $\Delta$ is homeomorphic to $\Zar(k(t)|k)$, when both sets are endowed with the constructible topology. Hence, it is enough to prove that $\Zar(k(t)|k)^\cons$ is not metrizable.

The points of $\Zar(k(t)|k)$ are $k(t)$, $k[t^{-1}]_{(t^{-1})}$ and the valuation rings of the form $k[t]_{(f(t))}$, where $f\in k[t]$ is an irreducible polynomial. The points different from $k(t)$ are isolated: indeed, $k[t^{-1}]_{(t^{-1})}$ is the only point in the open set $\Zar(k(t)|k)\setminus B(t)$, while $k[t]_{(f(t))}$ is the only point in the open set $\Zar(k(t)|k)\setminus B(f(t)^{-1})$. Since $\Zar(k(t)|k)^\cons$ is compact the claim follows from  Lemma \ref{lemma:metriz-uncount-limit}.
\end{proof}

In the following, we study more deeply spaces like $\{Z\in\Zar(L|V)\mid Z\subseteq W\}$ by using two classes of sequences that are similar to pseudo-convergent sequences. Let $E=\{s_n\}_{n\inN}$ be a sequence in $K$; then, we say that:
\begin{itemize}
\item $E$ is a \emph{pseudo-divergent sequence} if $v(s_n-s_{n+1})>v(s_{n+1}-s_{n+2})$ for every $n\inN$;
\item $E$ is a \emph{pseudo-stationary sequence} if $v(s_n-s_m)=v(s_{n'}-s_{m'})$ for every $n\neq m$, $n'\neq m'$.
\end{itemize}
These two kinds of sequences have been introduced in \cite{ChabPolCloVal} and together with the class of pseudo-convergent sequences introduced by Ostrowski form the class of \emph{pseudo-monotone sequences} \cite{PerPrufer}. Most of the notions introduced for pseudo-convergent sequences, like the breadth and the valuation domain $V_E$, can be generalized to pseudo-monotone sequences, see \cite{PS2}. In particular, the notion of pseudo-limit generalizes as well; however, there are fewer subsets that can be the set $\limiti_E$ of pseudo-limits of $E$. More precisely:
\begin{itemize}
\item if $E$ is a pseudo-divergent sequence, then there is an $\alpha\in K$ such that $\limiti_E=\{x\in K\mid v(x-\alpha)>\delta_E\}$, where $\delta_E$ is the breadth of $E$; if $\delta_E=v(c)\in\Gamma_v$, in particular, $\limiti_E=\alpha+cM$;
\item if $E$ is pseudo-stationary sequence, then there is an $\alpha\in K$ such that $\limiti_E=\alpha+cV$, where $v(c)=\delta_E$ is the breadth of $E$.
\end{itemize}
Furthermore, for every set $\limiti$ of this kind (with the additional hypothesis that the residue field of $V$ is infinite for pseudo-stationary sequences) there is a sequence $E$ of the right type with $\limiti=\limiti_E$. In particular, both pseudo-divergent sequences and pseudo-stationary sequences always have a pseudo-limit in $E$, and the elements of $E$ themselves are pseudo-limits of $E$ (\cite[Lemma 2.5]{PS2}). For both pseudo-divergent and pseudo-stationary sequences the ring $V_E$ is uniquely determined by the pseudo-limits: i.e., if $E,F$ are pseudo-divergent (respectively, pseudo-stationary) then $V_E=V_F$ if and only if $\limiti_E=\limiti_F$. For details, see \cite[Section 2.4]{PS2}.

Suppose that the residue field $k$ of $V$ is infinite, and let $Z=\{z_t\}_{t\in k}$ be a complete set of residues of $k$. Fix two elements $\alpha,c\in K$, and let $\delta=v(c)$. Let $\limiti=\{x\in K\mid v(x-\alpha)\geq\delta\}=\alpha+cV$ be the closed ball of center $\alpha$ and radius $\delta$. Then, there are a pseudo-convergent sequence $E$ and a pseudo-stationary sequence $F$ such that $\limiti_E=\limiti=\limiti_F$; by \cite[Proposition 7.1]{PS2}, $V_E\subsetneq V_F$. 

For every $z\in Z$, there is also a pseudo-divergent sequence $D_z$ such that $\limiti_{D_z}=\alpha-cz+cM$. Then, $V_{D_z}\neq V_E$ and $V_{D_z}\subsetneq V_F$ for every $z\in Z$; furthermore, $V_{D_z}\neq V_{D_{z'}}$ if $z\neq z'$. Let
\begin{equation*}
\mathcal{X}_{\alpha,\delta}=\{V_E,V_F,V_{D_z}\mid z\in Z\}
\end{equation*}
be the set of the rings in this form. By \cite[Proposition 7.2]{PS2}, the map $\widetilde{\pi}$ of Lemma \ref{lemma:quoziente} restricts to
\begin{equation*}
\begin{aligned}
\widetilde{\pi}\colon\mathcal{X}_{\alpha,\delta} & \longrightarrow\Zar(k(t)|k)\\
V_F & \longmapsto k(t),\\
V_E & \longmapsto k[1/t]_{(1/t)},\\
V_{D_z} & \longmapsto k[t]_{(t-\pi(z))},
\end{aligned}
\end{equation*}
and the lemma guarantees that $\widetilde{\pi}$ is also a homeomorphism between $\mathcal{X}_{\alpha,\delta}$ and its image.

In particular, we get the following; we denote by $\insiemeVEdiv$ the set of valuation rings $V_E$, as $E$ ranges among the pseudo-divergent sequences.

\begin{prop}\label{prop:VEdiv-delta}
Let $\insiemeVEdiv(\bullet,\delta)=\{V_E\mid E$ is a pseudo-divergent sequence with $\delta_E=\delta\}$. Then:
\begin{enumerate}[(a)]
\item\label{prop:VEdiv-delta:notGammav} if $\delta\notin\Gamma_v$, then $\insiemeVEdiv(\bullet,\delta)=\insiemeVE_K(\bullet,\delta)$;
\item\label{prop:VEdiv-delta:resfinito} if $\delta\in\Gamma_v$ and the residue field of $V$ is finite, then $\insiemeVEdiv(\bullet,\delta)$ is discrete (with respect to the Zariski and the constructible topology);
\item\label{prop:VEdiv-delta:resifinito} if $\delta\in\Gamma_v$ and the residue field of $V$ is infinite, then $\insiemeVEdiv(\bullet,\delta)$ is not Hausdorff (with respect to the Zariski topology).
\end{enumerate}
In particular, if the residue field of $V$ is infinite then $\insiemeVEdiv$ is not Hausdorff, with respect to the Zariski topology.
\end{prop}
\begin{proof}
\ref{prop:VEdiv-delta:notGammav} If $\delta\notin\Gamma_v$, then for every $\beta\in K$ we have $\{x\in K\mid v(x-\beta)\geq\delta\}=\{x\in K\mid v(x-\beta)>\delta\}$. Hence, if $E,F$ are, respectively, a pseudo-convergent and a pseudo-divergent sequence having $\beta$ as a pseudo-limit and having breadth $\delta$ then $\limiti_E=\limiti_F$, and thus by \cite[Proposition 5.1]{PS2} $V_E=V_F$. Since every pseudo-divergent sequence has pseudo-limits in $K$, it follows that $\insiemeVEdiv(\bullet,\delta)=\insiemeVE_K(\bullet,\delta)$.

\ref{prop:VEdiv-delta:resfinito} Suppose that $\delta\in\Gamma_v$, and let $c\in K$ be such that $v(c)=\delta$. Let $E$ be a pseudo-divergent sequence with breadth $\delta$, and let $\alpha\in\limiti_E$; since the residue field is finite we can find $\beta_1,\ldots,\beta_k\in K$ such that $0,\frac{\alpha-\beta_1}{c},\ldots,\frac{\alpha-\beta_k}{c}$ is a complete set of residues of the residue field of $V$.

We claim that
\begin{equation*}
\{V_E\}=B\left(\frac{X-\alpha}{c}\right)\cap B\left(\frac{c}{X-\beta_k}\right)\cap\cdots\cap B\left(\frac{c}{X-\beta_1}\right)\cap\insiemeVEdiv(\bullet,\delta).
\end{equation*}
Let $\Omega$ be the intersection on the right hand side. Since $\alpha\in\limiti_E$ the value $v(s_n-\alpha)$ decreases to $\delta$, and thus $v\left(\frac{s_n-\alpha}{c}\right)$ is always positive; in particular, $V_E\in B\left(\frac{X-\alpha}{c}\right)$. On the other hand, $v(s_n-\beta_i)=v(s_n-\alpha+\alpha-\beta_i)=v(\alpha-\beta_i)=\delta$ for every $i\in\{1,\ldots,k\}$ and every $n$, and thus $v\left(\frac{c}{s_n-\beta_i}\right)=0$, i.e., $V_E\in B\left(\frac{c}{X-\beta_1}\right)$. Hence, $V_E\in\Omega$.

Suppose now that $F=\{t_n\}_{n\inN}$ is a pseudo-divergent sequence such that $V_F\in\Omega$. Then, $V_F\in B\left(\frac{X-\alpha}{c}\right)$, i.e., $v(t_n-\alpha)\geq\delta$ for all large $n$, and thus $F$ must be eventually contained in the closed ball $\{x\in K\mid v(x-\alpha)\geq\delta\}=\alpha+cV=\beta_i+cV$ (for every $i$). Since $F$ has breadth $\delta$, by the discussion after Proposition \ref{prop:uncount-resfield} its set $\limiti_F$ of pseudo-limits is in the form $z+cM_V$, where $z$ is any element of $\limiti_F$; therefore, $\limiti_F$ is either $\alpha+cM_V$ or $\beta_i+cM_V$ for some $i$ (by the assumption on the $\beta_i$'s). However, if $\limiti_F=\beta_i+cM_V$ then $v(t_n-\beta_i)>\delta$ for all $n$, which implies that $V_F\notin B\left(\frac{c}{X-\beta_i}\right)$, against $V_F\in\Omega$; therefore, $\limiti_F=\alpha+cM_V=\limiti_E$ and thus $V_F=V_E$ by by \cite[Proposition 5.1]{PS2}. Therefore, $\Omega=\{V_E\}$ and $V_E$ is isolated. Since $V_E$ was arbitrary, $\insiemeVEdiv(\bullet,\delta)$ is discrete.

\ref{prop:VEdiv-delta:resifinito} Suppose that $\delta\in\Gamma_v$ and that the residue field is infinite. With the notation as before the statement, consider the set $\mathcal{X}_d(\alpha,\delta)=\mathcal{X}_{\alpha,\delta}\setminus\{V_F,V_E\}$: then, $\mathcal{X}_d(\alpha,\delta)$ is a subset of $\insiemeVEdiv(\bullet,\delta)$, and by Lemma \ref{lemma:quoziente} it is homeomorphic to $\Lambda=\{k[t]_{(t-z)}\mid z\in k\}\subseteq\Zar(k(t)|k)$. The map $\Spec(k[t])\longrightarrow\Zar(k[t])$, $P\mapsto K[t]_P$, is a homeomorphism (when both sets are endowed with the respective Zariski topologies) \cite[Lemma 2.4]{dobbs-fedder-fontana}, and thus $\Lambda$ is homeomorphic to $\Lambda_s=\{(t-z)\mid z\in k\}\subseteq\Max(k[t])$. The Zariski topology on $\Max(k[t])$ coincides with the cofinite topology; since $\Lambda_s$ is infinite, it follows that $\Lambda_s$ is not Hausdorff; thus, neither $\Lambda$ nor $\mathcal{X}_d(\alpha,\delta)$ nor $\insiemeVEdiv(\bullet,\delta)$ are Hausdorff.
\end{proof}

On the other hand, if we fix a pseudo-limit, we obtain a situation very similar to the pseudo-convergent case.
\begin{prop}\label{prop:VEdiv-beta}
Let $\beta\in K$, and let $\insiemeVEdiv(\beta,\bullet)=\{V_E\mid E$ is a pseudo-divergent sequence with $\beta\in\limiti_E\}$. Then,
\begin{equation*}
\insiemeVEdiv(\beta,\bullet)\simeq\insiemeVE(\beta,\bullet)\simeq(-\infty,+\infty]_{\insQ\Gamma_v}.
\end{equation*}
\end{prop}
\begin{proof}
For every $\beta,\beta'\in K$, we have $\insiemeVEdiv(\beta,\bullet)\simeq\insiemeVEdiv(\beta',\bullet)$ and $\insiemeVE(\beta,\bullet)\simeq\insiemeVE(\beta',\bullet)$, so we can suppose $\beta=0$.

Consider the map
\begin{equation*}
\begin{aligned}
\psi\colon K(X) & \longrightarrow K(X)\\
\phi(X) & \longmapsto \phi(1/X).
\end{aligned}
\end{equation*}

Then, $\psi$ is a $K$-automorphism of $K(X)$ that coincide with its own inverse, and thus it induces a self-homeomorphism
\begin{equation*}
\begin{aligned}
\overline{\psi}\colon \Zar(K(X)|V) & \longrightarrow \Zar(K(X)|V)\\
V_E & \longmapsto \psi(V_E).
\end{aligned}
\end{equation*}
We claim that $\overline{\psi}$ sends $\insiemeVEdiv(0,\bullet)$ to $\insiemeVE(0,\bullet)$, and conversely.

Note first that, for every $\phi\in K(X)$ and every $t\in K$, we have $\phi(t)=(\psi(\phi))(t^{-1})$.

Suppose $E=\{s_n\}_{n\inN}$ is a pseudo-divergent sequence having $0$ as a pseudo-limit; without loss of generality, $0\neq s_n$ for every $n$. Then, $\delta_n=v(s_n)$ is decreasing, and thus $F=\{s_n^{-1}\}_{n\inN}$ is a pseudo-convergent sequence having $0$ as a pseudo-limit. Then, $\phi(s_n)=(\psi(\phi))(s_n^{-1})$, and thus $\phi\in V_E$ if and only if $\psi(\phi)\in V_F$, i.e., $\psi(V_E)=V_F$, so that $\overline{\psi}(\insiemeVEdiv(0,\bullet))\subseteq\insiemeVE(0,\bullet)$. Conversely, if $F=\{t_n\}_{n\inN}$ is a pseudo-convergent sequence having $0$ as a pseudo-limit, then $E=\{t_n^{-1}\}_{n\inN}$ is a pseudo-divergent sequence with $0\in\limiti_E$, and as above $\phi\in V_F$ if and only if $\psi(\phi)\in V_E$, i.e., $\overline{\psi}(\insiemeVE(0,\bullet))\subseteq\insiemeVEdiv(0,\bullet)$. 

Since $\overline{\psi}$ is idempotent, it follows that $\overline{\psi}(\insiemeVE(0,\bullet))=\insiemeVEdiv(0,\bullet)$, and so $\insiemeVEdiv(0,\bullet)$ and $\insiemeVE(0,\bullet)$ are homeomorphic. The homeomorphism $\insiemeVE(0,\bullet)\simeq(-\infty,+\infty]_{\insQ\Gamma_v}$ follows from Theorem \ref{teor:betafix}.
\end{proof}

Note that, while the homeomorphism between $\insiemeVE(\beta,\bullet)$ and $(-\infty,+\infty]_{\insQ\Gamma_v}$ is constructed by sending $V_E$ to $\delta_E$ (Theorem \ref{teor:betafix}), the one between $\insiemeVEdiv(\beta,\bullet)$ and $(-\infty,+\infty]_{\insQ\Gamma_v}$ sends $V_E$ to $-\delta_E$.

We conclude with analyzing the pseudo-stationary case, showing that the two partitions give rise to especially uninteresting spaces.
\begin{prop}\label{prop:VEstaz}
The following hold.
\begin{enumerate}[(a)]
\item\label{prop:VEstaz:delta} For every $\delta\in\Gamma_v$, the set
\begin{equation*}
\insiemeVEstaz(\bullet,\delta)=\{V_E\mid E\text{~is a pseudo-stationary sequence with~}\delta_E=\delta\}
\end{equation*}
is discrete, with respect to the Zariski and the constructible topology.
\item\label{prop:VEstaz:beta} For every $\beta\in K$, the set
\begin{equation*}
\insiemeVEstaz(\beta,\bullet)=\{V_E\mid E\text{~is a pseudo-stationary sequence with~}\beta\in\limiti_E\}
\end{equation*}
is discrete, with respect to the Zariski and the constructible topology.
\end{enumerate}
\end{prop}
\begin{proof}
Since the constructible topology is finer than the Zariski topology, it is enough to prove the claim for the latter.

\ref{prop:VEstaz:delta} Take a pseudo-stationary sequence $E=\{s_n\}_{n\inN}$ of breadth $\delta$, and let $\beta\in\limiti_E$; let also $c\in K$ be such that $v(c)=\delta$. Consider the function $\phi(X)=\frac{X-\beta}{c}$; we claim that $B(\phi)\cap\insiemeVEstaz(\bullet,\delta)=\{V_E\}$.

Indeed, for large $n$ we have $v(s_n-\beta)=\delta$, and thus $v(\phi(s_n))=v(s_n-\beta)-v(c)=0$, so that $\phi\in V_E$, i.e., $V_E\in B(\phi)$. Conversely, suppose $V_F\in B(\phi)$, where $F=\{t_n\}_{n\inN}$ is pseudo-stationary with breadth $\delta$. Then, for large $n$, we must have $v(t_n-\beta)\geq\delta$. Since $v(t_n-t_m)=\delta$ for $n\neq m$, we must have $v(t_n-\beta)=\delta$, i.e., $\beta$ is a pseudo-limit of $F$. Thus, $\limiti_E=\beta+cV=\limiti_F$ and $V_E=V_F$ by \cite[Proposition 5.1]{PS2},

Therefore, $B(\phi)\cap\insiemeVEstaz(\bullet,\delta)=\{V_E\}$ and $V_E$ is an isolated point of $\insiemeVEstaz(\bullet,\delta)$. Since $V_E$ was arbitrary, $\insiemeVEstaz(\bullet,\delta)$ is discrete, as claimed.

\ref{prop:VEstaz:beta} Let $E=\{s_n\}_{n\inN}$ be a pseudo-stationary sequence having $\beta$ as a pseudo-limit, and let $c\in K$ be such that $v(c)=\delta_E$. Let $\phi(X)=\frac{X-\beta}{c}$; we claim that $B(\phi,\phi^{-1})\cap\insiemeVEstaz(\beta,\bullet)=\{V_E\}$.

The proof that $V_E\in B(\phi,\phi^{-1})$ follows as in the previous case. Suppose now that $F=\{t_n\}_{n\inN}$ is in the intersection. Then, we must have $v(\phi(t_n))\geq 0$ and $v(\phi^{-1}(t_n))=-v(\phi(t_n))\geq 0$; thus, $v(t_n-\beta)=\delta_E$ for large $n$. However, since $\beta$ is a pseudo-limit of $F$, we also have $v(t_n-\beta)=\delta_F$; hence, $\delta_E=\delta_F$ and $V_E=V_F$. Therefore, as above, $V_E$ is an isolated point of $\insiemeVEstaz(\beta,\bullet)$, which thus is discrete.
\end{proof}


\begin{thebibliography}{99}

\bibitem{APZ1} V. Alexandru, N. Popescu, Al. Zaharescu, \emph{A theorem of characterization of residual transcendental extensions of a valuation}. J. Math. Kyoto Univ. 28 (1988), no. 4, 579-592.

\bibitem{APZ2} V. Alexandru, N. Popescu, Al. Zaharescu, \emph{Minimal pairs of definition of a residual transcendental extension of a valuation}, J. Math. Kyoto Univ. 30 (1990), no. 2, 207-225. 

\bibitem{Bourb} N. Bourbaki, {\em Alg\`ebre commutative}, Hermann, Paris, 1961.


\bibitem{ChabPolCloVal} J.-L. Chabert, \emph{On the polynomial closure in a valued field}, J. Number Theory 130 (2010), 458-468.

\bibitem{dobbs-fedder-fontana} D. E. Dobbs, R. Fedder and M. Fontana, \emph{Abstract Riemann surfaces of integral domains and spectral spaces}, Ann. Mat. Pura Appl. (4) \textbf{148} (1987) 101--115.

\bibitem{gilmer} R. Gilmer, \emph{Multiplicative ideal theory}, Marcel Dekker Inc., New York, 1972, Pure and Applied Mathematics, No. 12.


\bibitem{hochster_spectral}
M. Hochster, \emph{Prime ideal structure in commutative rings}, Trans. Amer. Math. Soc., 142 (1969), 43--60.

\bibitem{Kap} I. Kaplansky, \emph{Maximal Fields with Valuation}, Duke Math. J. 9 (1942), 303-321.


\bibitem{maclane} S. MacLane, \emph{A construction for absolute values in polynomial rings}, Trans. Amer. Math. Soc. 40 (1936), no. 3, 363–395.

\bibitem{Ostr} A. Ostrowski, \emph{Untersuchungen zur arithmetischen Theorie der K\"orper},
Math. Z. 39 (1935), 269-404.

\bibitem{PerPrufer} G. Peruginelli, \emph{Pr\"ufer intersection of valuation domains of a field of rational functions}, J. Algebra 509 (2018), 240-262.

\bibitem{PerTransc} G. Peruginelli, \emph{Transcendental extensions of a valuation domain of rank one}, Proc. Amer. Math. Soc. 145 (2017), no. 10, 4211-4226.

\bibitem{PS1} G. Peruginelli, D. Spirito, \emph{The Zariski-Riemann space of valuation domains associated to pseudo-convergent sequences}, Trans. Amer. Math. Soc. 373 (2020), no. 11, 7959–7990.


\bibitem{PS2} G. Peruginelli, D. Spirito, \emph{Extending valuations to the field of rational functions using pseudo-monotone sequences}, J. Algebra 586 (2021), 756-786.



\bibitem{isolated-Zariski} D. Spirito, \emph{Isolated points of the Zariski space}, preprint,  arXiv: \href{https://arxiv.org/abs/2009.11141}{https://arxiv.org/abs/2009.11141}.

\bibitem{counterexamples-topology} L. A. Steen and J. A. Seebach, \emph{Counterexamples in Topology}. Springer-Verlag, New York-Heidelberg, second edition, 1978.

\bibitem{vaquie} M. Vaquié, Michel, \emph{Extension d'une valuation}, Trans. Amer. Math. Soc. 359 (2007), no. 7, 3439–3481. 

\bibitem{Warner} S. Warner, \emph{Topological Fields}. North-Holland Mathematics Studies, 157.

\bibitem{zariski_sing} O. Zariski, \emph{The reduction of the singularities of an algebraic surface}, Ann. of Math. (2), 40:639--689, 1939.

\bibitem{zariski_comp} O. Zariski, \emph{The compactness of the {R}iemann manifold of an abstract field of algebraic functions}, Bull. Amer. Math. Soc., 50:683--691, 1944.

\bibitem{ZS2} O. Zariski, P. Samuel, \emph{Commutative Algebra, vol. II}, Springer-Verlag, New York-Heidelberg-Berlin, 1975.

\end{thebibliography}
\end{document}